\documentclass[12pt]{article}
\usepackage{amssymb,amsmath} % font used for R in Real numbers
\newtheorem{theorem}{Theorem}[section]

\newtheorem{lemma}[theorem]{Lemma}

\numberwithin{equation}{section}
\newcommand \zm{\noindent{\bf Proof. \ }}
\newcommand \ezm{\hfill $\square$\par}

\title{Existence results for a superlinear singular equation of Caffarelli-Kohn-Nirenberg type}

\author{{Benjin Xuan}
\thanks{Supported by Grant 10101024 and 10371116 from the National Natural Science Foundation of China.}\\
{\it Department of Mathematics}\\
{\it University of Science and Technology of China}\\
{\it Universidad Nacional de Colombia}\\
{\it e-mail:wenyuanxbj@yahoo.com}}

\date{}
\begin{document}
\maketitle

\begin{abstract}
In this paper, using Mountain Pass Lemma and Linking Argument, we
prove the existence of nontrivial weak solutions for the Dirichlet
problem for the superlinear equation of Caffarelli-Kohn-Nirenberg
type in the case where the parameter $\lambda\in (0, \lambda_2)$,
$\lambda_2$ is the second positive eigenvalue of the quasilinear
elliptic equation of Caffarelli-Kohn-Nirenberg type.

\noindent{\bf Key Words:} singular equation, Caffarelli-Kohn-Nirenberg inequality, Mountain Pass Lemma, Linking Argument

\noindent{\bf Mathematics Subject Classifications:} 35J60.
\end{abstract}
\section{Introduction.}

In this paper, we shall investigate the existence of weak solutions for the following Dirichlet problem for the superlinear singular equation of Caffarelli-Kohn-Nirenberg type:
\begin{equation}
\label{eq1.1}
\left\{
\begin{array}{l}
-\mbox{div\,}(|x|^{-ap}|Du|^{p-2}Du)=\lambda |x|^{-(a+1)p+c}|u|^{p-2}u+|x|^{-bq}f(u),\mbox{ in } \Omega\\[2mm]
u=0,\  \mbox{on\,} \partial \Omega,
\end{array}
\right.
\end{equation}
where $\Omega\subset \mathbb{R}^n$ is an open bounded domain with $C^1$ boundary and $0\in \Omega$, $1<p<n,\ 0\leq a<\frac{n-p}p,\ a\leq b\leq a+1,\ q<p^*(a,b)=\frac{np}{n-dp},\ d=1+a-b\in [0,\ 1],\ c>0$.

For $a=0,\ c=p$, many results on linking-type critical points have been obtained (eg. \cite{AR, BR, CK} for $p=2$, \cite{FL} for $p \neq 2$ and \cite{XB} for the case with indefinite weights).

The starting point of the variational approach to these problems with $a\geq 0$ is the following weighted Sobolev-Hardy inequality due to Caffarelli, Kohn and Nirenberg \cite{CKN}, which is called the Caffarelli-Kohn-Nirenberg inequality. Let $1<p<n$. For all $u\in C_0^\infty(\mathbb{R}^n)$, there is a constant $C_{a,b}>0$ such that
\begin{equation}
\label{eq1.2}
\Big(\int_{\mathbb{R}^n}|x|^{-bq}|u|^{q}\,dx \Big)^{p/q}\leq C_{a,b}\int_{\mathbb{R}^n}|x|^{-ap}|Du|^{p}\,dx,
\end{equation}
where
\begin{equation}
\label{eq1.3}
-\infty< a<\frac{n-p}p,\ a\leq b\leq a+1,\ q=p^*(a,b)=\frac{np}{n-dp},\ d=1+a-b.
\end{equation}

Let $\Omega\subset \mathbb{R}^n$ is an open bounded domain with $C^1$ boundary and $0\in \Omega$, ${\cal D}_a^{1,p}(\Omega)$ be the completion of $C_0^\infty(\mathbb{R}^n)$, with respect to the norm $\|\cdot\|$ defined by
$$
\|u\|=\Big(\int_{\Omega}|x|^{-ap}|Du|^{p}\,dx \Big)^{1/p}.
$$
From the boundedness of $\Omega$ and the standard approximation
argument, it is easy to see that (\ref{eq1.2}) holds for any $u\in
{\cal D}_a^{1,p}(\Omega)$ in the sense:
\begin{equation}
\label{eq1.4}
\Big(\int_{\Omega}|x|^{-\alpha}|u|^{r}\,dx \Big)^{p/r}\leq C \int_{\Omega}|x|^{-ap}|Du|^{p}\,dx,
\end{equation}
for $1\leq r\leq \frac{np}{n-p},\ \alpha \leq (1+a)r+n(1-\frac rp)$, that is, the embedding ${\cal D}_a^{1,p}(\Omega) \hookrightarrow L^r(\Omega, |x|^{-\alpha})$ is continuous, where $L^r(\Omega, |x|^{-\alpha})$ is the weighted $L^r$ space with norm:
$$
\|u\|_{r, \alpha}:=\|u\|_{L^r(\Omega, |x|^{-\alpha})}=\Big( \int_{\Omega}|x|^{-\alpha}|u|^{r}\,dx\Big)^{1/r}.
$$

In fact, we have the following compact embedding result which is an extension of the classical Rellich-Kondrachov compactness theorem (cf. \cite{CC} for $p=2$ and  \cite{XB2} for the general case). For the convenience of readers, we include the proof here.
\begin{theorem}[Compact embedding theorem]
\label{thm1.1}
Suppose that $\Omega\subset \mathbb{R}^n$ is an open bounded domain with $C^1$ boundary and $0\in \Omega$, $1<p<n,\ -\infty< a<\frac{n-p}p$. The embedding ${\cal D}_a^{1,p}(\Omega) \hookrightarrow L^r(\Omega, |x|^{-\alpha})$ is compact if $1\leq r< \frac{np}{n-p},\ \alpha < (1+a)r+n(1-\frac rp)$.
\end{theorem}
\zm The continuity of the embedding is a direct consequence of
the Caffarelli-Kohn-Nirenberg inequality (\ref{eq1.2}) or
(\ref{eq1.4}). To prove the compactness, let $\{u_m\}$ be a
bounded sequence in ${\cal D}_a^{1,p}(\Omega)$. For any $\rho>0$
with $B_\rho(0)\subset \Omega$ is a ball centered at the origin
with radius $\rho$, there holds $\{u_m\}\subset
W^{1,p}(\Omega\setminus B_\rho(0))$. Then the classical
Rellich-Kondrachov compactness theorem guarantees the existence of
a convergent subsequence of $\{u_m\}$ in $L^r(\Omega\setminus
B_\rho(0))$. By taking a diagonal sequence, we can assume without
loss of generality that $\{u_m\}$ converges in
$L^r(\Omega\setminus B_\rho(0))$ for any $\rho>0$.

On the other hand, for any $1\leq r< \frac{np}{n-p}$, there exists a $b\in (a, a+1]$ such that $r<q=p^*(a, b)=\frac{np}{n-dp},\ d=1+a-b\in [0,\ 1)$. From the Caffarelli-Kohn-Nirenberg inequality (\ref{eq1.2}) or (\ref{eq1.4}), $\{u_m\}$ is also bounded in $L^q(\Omega, |x|^{-bq})$. By the H\"{o}der inequality, for any $\delta>0$, there holds
\begin{equation}
\label{eq2.001}
\begin{array}{ll}
\displaystyle \int_{|x|<\delta}|x|^{-\alpha}|u_m-u_j|^{r}\,dx & \leq \Big( \displaystyle \int_{|x|<\delta}|x|^{-(\alpha-br)\frac q{q-r}}\,dx\Big)^{1-\frac rq}\\[3mm]
& \ \ \ \ \times
 \Big(\displaystyle \int_{\Omega}|x|^{-br}|u_m-u_j|^{r}\,dx\Big)^{r/q}\\[3mm]
& \leq C \Big( \displaystyle \int_0^\delta r^{n-1-(\alpha-br)\frac q{q-r}}\,dr\Big)^{1-\frac rq}\\[3mm]
& =C \delta ^{n-(\alpha-br)\frac q{q-r}},
\end{array}
\end{equation}
where $C>0$ is a constant independent of $m$. Since $\alpha
< (1+a)r+n(1-\frac rp)$, there holds $n-(\alpha-br)\frac
q{q-r}>0$. Therefore, for a given $\varepsilon>0$, we first fix
$\delta>0$ such that
$$
\int_{|x|<\delta}|x|^{-\alpha}|u_m-u_j|^{r}\,dx \leq \frac{\varepsilon}2, \ \forall\ m, j\in \mathbb{N}.
$$
Then we choose $N\in  \mathbb{N}$ such that
$$
\int_{\Omega\setminus B_\delta(0)}|x|^{-\alpha}|u_m-u_j|^{r}\,dx \leq C_\alpha \int_{\Omega\setminus B_\delta(0)} |u_m-u_j|^{r}\,dx \leq \frac{\varepsilon}2, \ \forall\ m, j\geq N,
$$
where $C_\alpha=\delta ^{-\alpha}$ if $\alpha\geq 0$ and $C_\alpha=(\mbox{diam\,}(\Omega) )^{-\alpha}$ if $\alpha< 0$. Thus
$$
\int_{\Omega}|x|^{-\alpha}|u_m-u_j|^{r}\,dx \leq \varepsilon , \ \forall\ m, j\geq N,
$$
that is, $\{u_m\}$ is a Cauchy sequence in $L^q(\Omega, |x|^{-bq})$.
\ezm

Our results will mainly rely on the results of the eigenvalue problem correspondent to problem (\ref{eq1.1}) in \cite{XB1}. Let us first recall the main results of \cite{XB1}. Consider the nonlinear eigenvalue problem:
\begin{equation}
\label{eq1.5}
\left\{
\begin{array}{l}
-\mbox{div\,}(|x|^{-ap}|Du|^{p-2}Du)=\lambda |x|^{-(a+1)p+c}|u|^{p-2}u,\mbox{ in } \Omega\\[2mm]
u= 0, \ \ \mbox{on } \partial\Omega,
\end{array}
\right.
\end{equation}
where $\Omega\subset \mathbb{R}^n$ is an open bounded domain with $C^1$ boundary and $0\in \Omega$, $1<p<n,\ 0\leq  a<\frac{n-p}p,\ c>0$.

Let us introduce the following functionals in ${\cal D}_a^{1,p}(\Omega)$:
$$
\Phi(u):= \int_{\Omega} |x|^{-ap}|Du|^{p}\,dx,\ \mbox{and } J(u):=\int_{\Omega} |x|^{-(a+1)p+c}|u|^p\,dx.
$$
For $c>0$, $J$ is well-defined. Furthermore, $\Phi, J\in C^1({\cal D}_a^{1,p}(\Omega),\mathbb{R})$, and a real value $\lambda$ is an eigenvalue of problem (\ref{eq1.5}) if and only if there exists $u\in {\cal D}_a^{1,p}(\Omega)\setminus\{0\}$ such that $\Phi^\prime(u)=\lambda J^\prime(u)$. At this point, let us introduce set
$$
{\cal M}:=\{u\in {\cal D}_a^{1,p}(\Omega)\ :\ J(u)=1 \}.
$$
Then ${\cal M}\neq \emptyset$ and ${\cal M}$ is a $C^1$ manifold in ${\cal D}_a^{1,p}(\Omega)$. It follows from the standard Lagrange multiples arguments that eigenvalues of (\ref{eq1.5}) correspond to critical values of $\Phi|_{{\cal M}}$. From Theorem \ref{thm1.1}, $\Phi$ satisfies the (PS) condition on ${\cal M}$. Thus a sequence of critical values of $\Phi|_{{\cal M}}$ comes from the Ljusternik-Schnirelman critical point theory on $C^1$ manifolds. Let $\gamma(A)$ denote the Krasnoselski's genus on ${\cal D}_a^{1,p}(\Omega)$ and for any $k\in \mathbb{N}$, set
$$
\Gamma_k:=\{A\subset {\cal M}\ :\ A \mbox{ is compact, symmetric and } \gamma(A)\geq k\}.
$$
Then values
\begin{equation}
\label{eq1.6}
\lambda_k:=\inf_{A\in\Gamma_k }\max_{u\in A} \Phi(u)
\end{equation}
are critical values and hence are eigenvalues of problem (\ref{eq1.5}). Moreover, $\lambda_1\leq \lambda_2\leq \cdots \leq   \lambda_k\leq \cdots \to+\infty$.

One can also define another sequence of critical values minimaxing $\Phi$ along a smaller family of symmetric subsets of ${\cal M}$. Let us denote by $S^k$ the unit sphere of $\mathbb{R}^{k+1}$ and
$$
{\cal O}(S^k, {\cal M}):= \{h\in C(S^k, {\cal M})\,:\, h \mbox{ is odd}\}.
$$
Then for any $k\in \mathbb{N}$, the value
\begin{equation}
\label{eq1.7}
\mu_k:= \inf_{h\in{\cal O}(S^{k-1}, {\cal M}) }\max_{t\in S^{k-1}} \Phi(h(t))
\end{equation}
is an eigenvalue of (\ref{eq1.5}). Moreover $\lambda_k\leq \mu_k$. This new sequence of eigenvalues was first introduced by \cite{DR} and later used in \cite{CM2, CM} for $a=0, c=p$.

In \cite{XB1}, we proved that the first positive eigenvalue $\lambda_1=\mu_1$ is simple, isolated and it is the unique eigenvalue with positive eigenfunction, and $\underline{\lambda}_2:=\inf\{\lambda\in \mathbb{R}\,:\, \lambda \mbox{ is eigenvalue and } \lambda>\lambda_1 \}=\lambda _2=\mu_2$.

In this paper, based on the Mountain Pass Lemma and the Linking Argument, we shall prove the existence of nontrivial weak solutions problem (\ref{eq1.1}) in the case where the parameter $\lambda\in (0, \lambda_2)$.

%%%%%%%%%%%%%%%%%%%%%%%%%%%%%%%%%%%%%%%%%%%%%%%%%%%%%%%%%%%%%%%%%%%%%%%%%%%%%%%%%%%%%%%%%%
\section{Linking results}

Let $e_k \in {\cal M}$ be the eigenfunction associated to $\lambda_k$, then $\| e_k\|_{{\cal D}_a^{1,p}(\Omega)}^p=\lambda_k$.
Denote $G=\{ u \in  {\cal M} :\ \Phi(u)< \lambda_2\}$.
Obviously, $G$ is an open set containing $e_1$ and $e_2$. Moreover $-G=G$. First we shall prove the following Lemma.

\begin{lemma}
\label{2.1}
$e_1$ and $-e_1$ do not belong to the same connected component of $G$.
\end{lemma}

{\bf Proof.} Otherwise, there exists a continuous curve $\sigma$ connecting $e_1$ and $-e_1$ in $G$. Let $A=\sigma \cup \{- \sigma\}$, then from the definition of ${\cal M}$, $0\not\in A$,
hence $\gamma(A)>1$. By connectedness of $A$, so $A \in \Gamma_2$. Hence, as $A$ is a compact set in $G$, and from the definition of $G$, we will have $\max\{\Phi(u);\ u \in A\}< \lambda_2$, and this contradicts the definition of $\lambda_2$.
 \hfill Q.E.D. \bigskip

Let $G_1$ be the connected component of $G$ containing $e_1$, then $-G_1$ is the connected component of $G$ containing $-e_1$. Let
$$
K_1=\{tu:\ u \in G_1, t>0 \},\ \  K=K_1 \cup \{-K_1\}.
$$
Then, we have
\begin{equation}
\label{eq2.1}
\int_{\Omega}|x|^{-ap}|D u|^p\,dx < \lambda_2 \int_{\Omega}|x|^{-(a+1)p+c} | u|^p\,dx, \ \forall u \in K,
\end{equation}
and
\begin{equation}
\label{eq2.2}
\int_{\Omega}|x|^{-ap}|D u|^p\,dx= \lambda_2 \int_{\Omega}|x|^{-(a+1)p+c} | u|^p\,dx, \ \forall u \in \partial K,
\end{equation}
where $\partial K$ is the boundary of $K$ in $X={\cal D}_a^{1,p}(\Omega)$. Let $(\partial K)_\rho=\{u \in \partial K:\ \|u\|=\rho\}$.

Set
$${\cal E}_1= \mbox{span}\, \{e_1\}, \ {\cal E}_2= \mbox{span}\, \{e_1,\ e_2\},$$
$$
{\cal Z}=\{u \in  X:\ \int_\Omega |D u|^p= \lambda_2 \int_\Omega V(x) |u|^p\}.
$$
(\ref{eq2.2}) implies $\partial K \subset {\cal Z}$.

Similar to Proposition 2.1-2.2 in \cite{FL} and Lemma 2.1-2.2 in \cite{XB}, we obtain the following two linking results.

\begin{theorem}
\label{2.2}
Assume that $v \in {\cal E}_1, \ v \neq 0,\ Q=[-v,\ v]$ is the line segment connecting $-v$ and $v$,
$\partial Q=\{-v,\ v\}$. Then $\partial Q \subset Q$ and ${\cal Z}$ link in $X$, that is,
\begin{enumerate}
\item[(i)] $\partial Q \cap {\cal Z} = \emptyset$ and
\item[(ii)] For any continuous map $\psi: Q \to X$ with $\psi |_{\partial Q}=\mbox{id}$, there holds $\psi(Q)\cap {\cal Z} \neq \emptyset$.
\end{enumerate}
\end{theorem}

\zm It is obvious that $\partial Q \cap {\cal Z} = \emptyset$. Now
let $\psi: Q=[-v, v]\to X$ be continuous and $\psi |_{\partial
Q}=\mbox{id}$. From the definition of $K$ and Lemma \ref{2.1}, $K$
has two connected components $K_1$ and $-K_1$. Assume $v\in K_1, \
-v\in -K_1$, then $\psi(Q)$ is a continuous curve connecting $v$
and $-v$, therefore there holds $\psi(Q)\cap \partial K\neq
\emptyset$ and thus $\psi(Q)\cap {\cal Z} \neq \emptyset$. \ezm

\begin{theorem}
\label{2.3}
Assume $0<\rho<r< \infty$, let $\tilde e_1=e_1/\lambda_1^{1/p},\ \tilde e_2=e_2/\lambda_2^{1/p}$, and
$$
Q=Q_r=\{u=t_1\tilde e_1+t_2\tilde e_2:\ \|u\|\leq r, t_2\geq 0 \},
$$
$$
\partial Q=\partial Q_r=\{u=t_1\tilde e_1:\  |t_1|\leq r \}\cup \{u \in Q_r:\  \|u\|= r\},
$$
$$
Z_\rho =\{u \in {\cal Z}:\ \|u\|=\rho \}.
$$
Then $\partial Q_r \subset Q_r$ and $Z_\rho $ link in $X$.
\end{theorem}

\zm
$\partial Q_r \cap  Z_\rho = \emptyset$ is obvious. Let $\psi: Q_r\to X$ be continuous and $\psi |_{\partial Q_r}=\mbox{id}$. Denote $d_1=\mbox{dist\,}(\tilde e_1,\partial K)$ and define mapping $P: X\to {\cal E}_2$ as follows:
$$
P(u)=\begin{cases} & \big(\min\{\mbox{dist\,}(u,\partial K), rd_1\}\big)\tilde e_1+(\|u\|-\rho) \tilde e_2,\ \mbox{ if }u \not\in -K_1;\\[2mm]
& -\big(\min\{\mbox{dist\,}(u,\partial K), rd_1\}\big)\tilde e_1+(\|u\|-\rho) \tilde e_2,\ \mbox{ if }u \in -K_1.
\end{cases}
$$
It is easy to see that $P$ is continuous, and $P$ maps $v=r \tilde e_1$ to $v_1=Pv=rd_1\tilde e_1+(r-\rho)\tilde e_2$, the origin $0$ to $0_1=P0=-\rho \tilde e_2$, the line segment $[v,0]$ onto the line segment $[v_1, 0_1]$ homeomorphically; $-v=-r\tilde e_1$ to $v_2=P(-v)=-rd_1\tilde e_1+(r-\rho)\tilde e_2$, the line segment $[0, -v]$ onto a line segment $[0_1, v_2]$
homeomorphically; and the half circle $\{u\in \partial Q: \ \|u\|=r \}$ which is from $v$ to$-v$
in $\partial Q$ onto the line segment $[v_1, v_2]$, where $P(r \tilde e_2)=(r-\rho)\tilde e_2$.

Let $f=P\circ \psi: Q\to {\cal E}_2$. From the discussion above,
there holds $0\not\in f(\partial Q)$, and when $u$ turns a circuit
along $\partial Q$ anticlockwise, $f(u)$ also moves a circuit
around the original $0$ in ${\cal E}_2$ anticlockwise. Hence by a
degree argument, there holds $\mbox{deg\,}(f, Q,0)=1$. So there
exists some $u\in Q$ such that $f(u)=0$, i.e., $P(\psi(u))=0$,
which implies that $\psi(u)\in \partial K$; and
$\|\psi(u)\|=\rho$. Thus $\psi(u) \in (\partial K )_\rho$ and
$\psi(Q)\cap (\partial K)_\rho \neq \emptyset$. Since $(\partial
K)_\rho \subset Z_\rho$, hence $\psi(Q)\cap {\cal Z} \neq
\emptyset$ \ezm

%%%%%%%%%%%%%%%%%%%%%%%%%%%%%%%%%%%%%%%%%%%%%%%%%%%%%%%%%%%%%%%%%%%%%%%%%%%%%%%%%%%%%%%%%%
\section{Existence results for problem (\ref{eq1.1})}

In this section, we will give some conditions on $f(u)$ to guarantee the functional associated to problem (\ref{eq1.1}) satisfies the Palais-Smale condition ((PS) condition) for $\lambda\in (0, \lambda_2)$, the geometric assumptions of Mountain Pass Lemma (cf. Theorem 6.1 in Chapter 2 of \cite{SM}) in the case of $0<\lambda<\lambda_1$, and those of the linking theorem (cf. Theorem 8.4 in Chapter 2 of \cite{SM}) in the case of $\lambda_1\leq \lambda <\lambda_2$.

Assume $f:\ \mathbb{R} \to \mathbb{R}$ satisfies:
\begin{enumerate}
\item[(f$_1$)] (Subcritical growth) $|f(s)|\leq c_1|s|^{q-1}+ c_2, \ \forall s \in \mathbb{R}$, where $1< q < p^*(a,b)=\frac{Np}{N-dp}$;

\item[(f$_2$)]$f\in C(\mathbb{R}, \mathbb{R})$, $f(0)=0,\ uf(u)\geq 0,\, u\in \mathbb{R}$;

\item[(f$_3$)](Asymptotic property at infinity) $\exists\, \theta\in(p,\ p^*(a,b))$ and $M>0$ such that $0< \theta F(u)\leq uf(u)$ for $|u|\geq M$, where $F(s)=\displaystyle \int_0^s f(t)dt$;

\item[(f$_4$)](Asymptotic property at $u=0$) $\lim\limits_{s \to 0}f(s)/|s|^{p-1}=0$.
\end{enumerate}

Theorem (\ref{thm1.1}) and (f$_1$) imply that functional $I:\ X \to \mathbb{R}$:
$$
I(u)=\frac1p\int_\Omega  |x|^{-ap} |D u|^p\,dx-\frac{\lambda}{p}\int_\Omega  |x|^{-(a+1)p+c} |u|^p\,dx
-\int_\Omega |x|^{-bq} F(u)dx
$$
is well-defined and $I \in C^1(X; \mathbb{R})$, and the weak solutions of problem (\ref{eq1.1}) is equivalent to the critical points of $I$. (f$_2$) implies that $0$ is a trivial solution to problem (\ref{eq1.1}).

\begin{lemma}
\label{3.01}
If $f$ satisfies assumptions (f$_1$)-(f$_3$), then
$I$ satisfies the (PS) condition for $\lambda\in (0, \lambda_1)$.
\end{lemma}
\noindent\textbf{Proof. 1. }The boundedness of $(PS)$ sequence of $I$.

Suppose $\{u_m\}$ is a $(PS)$ sequence of $I$, that is, there exists $C>0$ such that $|I(u_m)|\leq C$ and $I^\prime (u_m) \to 0$ in $X^\prime$, the dual space of $X$, as $m\to \infty$. The properties of the first eigenvalue $\lambda_1$ imply that for any $u\in X$,there holds
$$
\lambda_1\int_\Omega |x|^{-(a+1)p+c}|u|^p\,dx \leq \int_\Omega |x|^{-ap}  |D u|^p\,dx.
$$
Let $c:=\sup\limits_m I(u_{m})$. Then by the above inequality and (f$_3$), as $m\to \infty$, there holds
$$
\begin{array}{ll}
& c-\displaystyle\frac 1\theta o(1)\|u_{m}\|=(\displaystyle\frac1p -\displaystyle\frac 1\theta)\displaystyle\int_\Omega |x|^{-ap}  |D u_{m}|^p\,dx\\[3mm]
& \ \ \ \ \  \ \ \ \ \
-\lambda (\displaystyle\frac1p -\displaystyle\frac 1\theta)\displaystyle\int_\Omega |x|^{-(a+1)p+c} |u_{m}|^p \,dx
+ \displaystyle\int_\Omega |x|^{-bq}(\displaystyle\frac 1\theta f(u_{m})u_{m} - F(u_{m}))\,dx \\[3mm]
& \ \ \ \ \
 \geq (\displaystyle\frac1p- \displaystyle\frac 1\theta) (1- \displaystyle\frac{\lambda}{\lambda_1})\displaystyle\int_\Omega  |x|^{-ap}  |D u_{m}|^p\,dx
\\[3mm]
&  \ \ \ \ \  \ \ \ \ \
+ \displaystyle\int_{\Omega(u_{m}\geq M)} |x|^{-bq} (\displaystyle\frac 1\theta f(u_{m})u_{m} - F(u_{m}))\,dx\\[3mm]
&  \ \ \ \ \  \ \ \ \ \
+ \displaystyle\int_{\Omega(u_{m}< M)}|x|^{-bq}(\displaystyle\frac 1\theta f(u_{m})u_{m} - F(u_{m}))\,dx \\[3mm]
& \ \ \ \ \
 \geq (\displaystyle\frac1p -\displaystyle\frac 1\theta) (1- \displaystyle\frac{\lambda}{\lambda_1}) \|u_{m}\|^p-C_1,
\end{array}
$$
where $C_1\geq 0$ is a constant independent of $u_{m}$. The above estimate
implies the boundedness of $\{u_{m}\}$ for $0<\lambda < \lambda_1$.

\noindent\textbf{2. }By (f$_1$), $f$ satisfies the subcritical growth condition, by a standard argument, one can obtain that there exists a convergent subsequence of $\{u_{m}\}$ from the boundedness of $\{u_{m}\}$ in $X$.
\ezm

\begin{theorem}
\label{3.1}
If $f$ satisfies assumptions (f$_1$)-(f$_4$), then problem (\ref{eq1.1}) has a nontrivial weak solution $u \in W_0^{1,p}(\Omega)$ provided that $0<\lambda <\lambda_1$.
\end{theorem}
\noindent\textbf{Proof. }We will verify the geometric assumptions of the Mountain Pass Lemma (cf. \cite{SM} Chapter 2, Theorem 6.1):

\begin{enumerate}
\item[(1)] $I(0)=0$ is obvious;
\item[(2)] $\exists \rho >0,\ \alpha >0:\ \|u\|=\rho \implies I(u)\geq \alpha$;
\end{enumerate}

In fact, $\forall\,u \in X$, there holds
\begin{equation}
\label{eq3.1}
I(u)\geq \displaystyle\frac1p(1-\frac{\lambda}{\lambda_1})\displaystyle\int_\Omega  |x|^{-ap}  |D u|^p\,dx-\displaystyle\int_\Omega |x|^{-bq} F(u)\,dx.
\end{equation}
From (f$_4$), $\forall\, \epsilon >0$, $\exists\, \rho_0=\rho_0(\epsilon)$ such that:
if $0< \rho=\|u\|< \rho_0$, then $|f(u)|< \epsilon |u|^{p-1}$, thus
$$
\int_\Omega |x|^{-bq} F(u)\,dx\leq \int_\Omega |x|^{-bq} \int_0^{u(x)}f(t)\,dt \,dx
\leq \frac{\epsilon}{p}\int_\Omega  |x|^{-bq}|u|^pdx\leq \frac{c_0\epsilon}{p} \|u\|.
$$
Choose $c_0\epsilon_0=(1- \frac{\lambda}{\lambda_1})/2 >0, \ \rho=\dfrac{\rho_0(\epsilon_0)}2$, from
(\ref{eq3.1}), one has
\begin{equation}
\label{eq3.2}
I(u)\geq \frac1p (1- \frac{\lambda}{\lambda_1} -c_0\epsilon_0)\int_\Omega |x|^{-ap}  |D u|^p\,dx
\geq \frac{\lambda_1-\lambda}{2\lambda_1 p}\cdot \rho=:\alpha >0.
\end{equation}

\begin{enumerate}
\item[(3)] $\exists\, u_1 \in X:\ \|u_1\|\geq \rho\ \mbox{and}\ I(u_1)< 0$.
\end{enumerate}

In fact, from (f$_2$) and (f$_3$), one can deduce that there exist constants $c_3, c_4>0$ such that
\begin{equation}
\label{eq3.03}
F(s)\geq c_3|s|^\theta -c_4, \ \forall s \in \mathbb{R}.
\end{equation}
Since $\theta>p $, a simple calculation shows that as $t\to \infty$, there holds
\begin{equation}
\label{eq3.04}
\begin{array}{ll}
I(te_1)&\leq \displaystyle\frac{t^p}p\displaystyle\int_\Omega |x|^{-ap}  |D e_1|^p\,dx-\displaystyle\frac{\lambda t^p}{p}\displaystyle\int_\Omega  |x|^{-(a+1)p+c} |e_1|^p\,dx\\[3mm]
& \ \ \ \ \ \ -c_3 t^\theta \displaystyle\int_\Omega |x|^{-bq}  |e_1|^\theta\,dx+c_4 \displaystyle\int_\Omega |x|^{-bq}\,dx\\[3mm]
&\to -\infty,
\end{array}
\end{equation}
which implies that $I(te_1)<0$ for $t>0$ large enough.

Thus Lemma \ref{3.01} and the Mountain Pass Lemma imply that value
$$
\beta=\inf_{p \in P}\sup_{u \in p}E(u)\geq \alpha>0
$$
is critical, where $P =\{p \in C^0([0,1];\, X):\ p(0)=0, p(1)=u_1\}$.
That is, there is a $u \in X$, such that
$$
E^\prime(u)=0,\ E(u)=\beta>0.
$$
$E(u)=\beta>0$ implies $u \not\equiv 0$.
\ezm

\begin{lemma}
\label{3.02}
Assume that $\lambda_1\leq \lambda< \lambda_2$ and $f$ satisfies assumptions (f$_1$)-(f$_3$). Then $I$ satisfies the (C)$_c$ condition introduced by Cerami in \cite{CG}, that is, any sequence $\{u_m\}\subset X$ such that $I(u_m)\to c$ and $(1+\|u_m\|)\|I^\prime (u_m)\|_{X^\prime}\to 0$ possesses a convergent subsequence.
\end{lemma}

\noindent\textbf{Proof. 1. }The boundedness of (C)$_c$ sequence in $X$.

Let $\{u_m\}\subset X$ be such that $I(u_m)\to c$ and $(1+\|u_m\|)\|I^\prime (u_m)\|_{X^\prime}\to 0$. Then from (f$_2$), (f$_3$) and (\ref{eq3.03}), as $m\to \infty$, there holds
\begin{equation}
\label{eq3.07}
\begin{array}{ll}
& pc+o(1)=pI(u_m)-<I^\prime (u_m), u_m>\\[3mm]
&\ \ \ \ =\displaystyle\int_\Omega|x|^{-bq} (u_mf(u_m)-pF(u_m))\,dx\\[3mm]
&\ \ \ \ =\displaystyle\int_\Omega |x|^{-bq} (u_mf(u_m)- \theta F(u_m))\,dx+(\theta-p)\displaystyle\int_\Omega \theta |x|^{-bq}F(u_m)\,dx\\[3mm]
&\ \ \ \ \geq -C_1+(\theta-p)c_3|u_m|_{L^\theta(\Omega, |x|^{-bq})}^\theta-C_4\displaystyle\int_\Omega |x|^{-bq}\,dx.
\end{array}
\end{equation}
Thus $\theta>p$ implies the boundedness of $\{u_m \}$ in $L^\theta(\Omega, |x|^{-bq})$.

On the other hand, a simple calculation shows that
\begin{equation}
\label{eq3.09}
\begin{array}{ll}
&\theta c+o(1)=\theta I(u_m)-<I^\prime (u_m), u_m>\\[3mm]
&\ \ \ \ =(\displaystyle\frac \theta p-1)\|D u_m\|_{L^p(\Omega, |x|^{-ap})}^p-\lambda (\displaystyle\frac \theta p-1)\displaystyle\int_\Omega |x|^{-(a+1)p+c} |u_m|^p\,dx\\[3mm]
&\ \ \ \ \ \ \ \ \ +   \displaystyle\int_\Omega  |x|^{-bq}(u_mf(u_m)-\theta F(u_m))\,dx\\[3mm]
&\ \ \ \ \geq (\displaystyle\frac \theta p-1)\displaystyle\int_\Omega |x|^{-ap}|D u_m|^p\,dx -C \\[3mm]
&\ \ \ \ \ \ \ \ \ +  \displaystyle\int_{\Omega(u_m< M)} |x|^{-bq} (u_mf(u_m)- \theta F(u_m))\,dx\\[3mm]
&\ \ \ \ \ \ \ \ \  + \displaystyle\int_{\Omega(u_m\geq M)}|x|^{-bq} (u_mf(u_m)- \theta F(u_m))\,dx\\[3mm]
&\ \ \ \ \geq (\displaystyle\frac \theta p-1) \|D u_m\|_{L^p(\Omega, |x|^{-ap})}^p -C,
\end{array}
\end{equation}
where $C>0$ is a universal constant independent of $u_m$, which may be different from line to line. Thus $\theta>p$ and (\ref{eq3.09}) imply the boundedness of $\{u_m \}$ in $X$.

\noindent\textbf{2. }By (f$_1$), $f$ satisfies the subcritical growth condition, by a standard argument, one can obtain that there exists a convergent subsequence of $\{u_{m}\}$ from the boundedness of $\{u_{m}\}$ in $X$.
\ezm

\begin{theorem}
\label{3.2}
Suppose that $f$ satisfies assumptions (f$_1$)-(f$_4$). Then problem (\ref{eq1.1}) has a nontrivial weak solution $u \in X$ provided that $\lambda_1\leq \lambda <\lambda_2$.
\end{theorem}
\noindent\textbf{Proof. }It was shown in \cite{BBF} that $(C)_c$ condition actually suffices to get a deformation theorem (Theorem 1.3 in \cite{BBF}, and it also remarked in \cite{CM1} that
the proofs of the standard Mountain Pass Lemma and saddle-point theorem go through without change once the deformation theorem (Theorem 1.3 in \cite{BBF} is obtained with $(C)_c$ condition. Here we verify the assumptions of standard Linking Argument Theorem(cf. \cite{SM} Chapter 2, Theorem 8.4) hold with $(C)_c$ condition replacing (PS)$_c$ condition.

Since $\partial Q_r \subset Q_r$ and $Z_\rho $ link in $X$, it suffice to show that:
\begin{enumerate}
\item[(1).]$\alpha_0=\sup\limits_{u\in \partial Q_r}I(u)\leq0$, when $r>0$ is large enough;
\item[(2).]$\alpha=\inf\limits_{u\in Z_\rho}I(u)>0$, when $\rho>0$ is small enough.
\end{enumerate}

In fact, let $u=te_1\in {\cal E}_1$, from assumption (f$_2$), $F(x,s)\geq 0$ for all $s\in \mathbb{R}$ and almost every $x\in \Omega$, thus there holds
\begin{equation}
\label{eq3.5}
\begin{array}{ll}
I(u)=I(te_1)&\leq \dfrac{|t|^p}{p}\displaystyle\int_\Omega |x|^{-ap} |D e_1|^p\,dx -\dfrac{|t|^p \lambda}{p}\displaystyle\int_\Omega |x|^{-(a+1)p+c}|e_1|^p\,dx\\[3mm]
& =\dfrac{|t|^p}{p}(1-\dfrac{\lambda}{\lambda_1})\|e_1\|\leq 0.
\end{array}
\end{equation}
Noticing that
$$
|u_m|_{L^\theta(\Omega, |x|^{-bq})}=\big(\int_\Omega |x|^{-bq} |u|^\theta\big)^{1/\theta},
$$
is a norm on ${\cal E}_2$, and the norms of finite dimensional space are equivalent, thus there exists a constant $c_5>0$ such that
$$
\int_\Omega |x|^{-bq}|u|^\theta \,dx \geq c_5 \|u\|^\theta ,
$$
From (\ref{eq3.03}), there holds
\begin{equation}
\label{eq3.6}
I(u)\leq \dfrac1p \|u\|^p-c_3c_5 \|u\|^\theta+c_4|\Omega|.
\end{equation}
Since $\theta >p$, there holds
$$
I(u)\to -\infty, \mbox{ as }\|u\|\to \infty, \ u\in {\cal E}_2.
$$
This implies (1).

From (f$_4$) and (f$_1$), there holds
$$
\int_\Omega  |x|^{-bq} F(u)\,dx=o(\|u \|^p) \mbox{ as } u\to 0 \mbox{ in }X,
$$
then for any $u\in Z$, there holds
\begin{equation}
\label{eq3.7}
I(u)= \frac1p(1-\frac{\lambda}{\lambda_2})\int_\Omega  |x|^{-ap}|D u|^p\,dx+o(\|u \|^p).
\end{equation}
Since $\lambda<\lambda_2$, (\ref{eq3.7}) implies (2).

Thus the Linking Argument Theorem (cf. \cite{SM} Chapter 2, Theorem 8.4) implies that value
$$
\beta=\inf_{h\in \Gamma}\sup_{u\in Q}E(h(u))\geq \alpha>0
$$
is critical, where $\Gamma =\{h\in C^0(X;X);\ h|_{\partial Q} =\mbox{id}\}$.
That is, there is a $u \in X$, such that
$$
E^\prime(u)=0,\ E(u)=\beta>0.
$$
$E(u)=\beta>0$ implies $u \not\equiv 0$.
\ezm

\noindent\textsc{Benjin Xuan} \\
Department of Mathematics\\
Universidad Nacional de Colombia, Bogot\'a, Colombia\\
University of Science and Technology of China, Hefei, P. R. China\\
e-mail: bjxuan@matematicas.unal.edu.co, wenyuanxbj@yahoo.com

\end{document}